\documentclass[10pt]{article}
\usepackage[tbtags]{amsmath}
\usepackage{epsfig,amstext,amssymb,amsthm,latexsym}
\allowdisplaybreaks[4]
\usepackage{amssymb,color}
\newtheorem{thm}{Theorem}
\newtheorem{cor}[thm]{Corollary}
\newtheorem{lem}[thm]{Lemma}

\theoremstyle{definition}
\newtheorem{defn}{Definition}
\newtheorem{exmp}{Example}
\newtheorem{remark}{Remark}

\newcommand{\BEX}{\begin{exmp}}
\newcommand{\EEX}{\end{exmp}}

\allowdisplaybreaks[4]
\newcommand{\nwc}{\newcommand}
\nwc{\COM}[1]{}

\newcommand{\nelem}[1]{{Lemma \ref{#1}}}

\newcommand{\netheo}[1]{{Theorem \ref{#1}}}

\newcommand{\kb}[1]{\boldsymbol{#1}}
\newcommand{\vk}[1]{\kb{#1}}

\def\FRE{\mbox{Fr\'{e}chet }}

\newcommand{\ve}{\varepsilon}
\newcommand{\abs}[1]{\lvert #1 \rvert}
\newcommand{\Abs}[1]{ \Bigl \lvert #1 \Bigr \rvert}

\newcommand{\E}[1]{\mbox{\rm$\vk{E}$}\{#1\}}
\newcommand{\pk}[1]{\mbox{\rm$\vk{P}$} \{#1\} }

\newcommand{\pb}[1]{\mbox{\rm$\vk{P}$}\Bigl \{#1 \Bigr \}}

\newcommand{\R}{\!I\!\!R}

\newcommand{\inr}{\in \R}

\newcommand{\ldot}{,\ldots,}

\newcommand{\limit}[1]{\lim_{#1 \to   \infty}}
\newcommand{\todis}{\stackrel{d}{\to}}
\newcommand{\toprob}{ \stackrel{p}{\to}}

\newcommand{\ntoi}{n \to \infty }
\newcommand{\equaldis}{\stackrel{d}{=}}

\newcommand{\BQN}{\begin{eqnarray}}
\newcommand{\EQN}{\end{eqnarray}}
\newcommand{\BQNY}{\begin{eqnarray*}}
\newcommand{\EQNY}{\end{eqnarray*}}
\newcommand{\BS}{\begin{sat}}
\newcommand{\ES}{\end{sat}}
\newcommand{\BRM}{\begin{remark}}
\newcommand{\ERM}{\end{remark}}

\newcommand{\BL}{\begin{lem}}
\newcommand{\EL}{\end{lem}}
\newcommand{\BT}{\begin{thm}}
\newcommand{\ET}{\end{thm}}
\newcommand{\BK}{\begin{cor}}
\newcommand{\EK}{\end{cor}}

\newcommand{\QED}{\hfill $\Box$}
\newcommand{\IF}{\infty}

\newcommand{\prooftheo}[1]{ \textsc{Proof of Theorem} \ref{#1} }

\def\wth{\widetilde h}
\def\btrho{\overline{h}_\rho}

\def\wthN{\widetilde h_n}
\def\btrhoN{\overline{h}_{n,\rho}}

\def\arhon{\alpha_n}
\def\aroxy{\alpha_\rho(x,y)}
\def\broxy{\alpha_\rho^*(x,y)}
\def\gnt{g_n(x)}

\def\wtiro{\widetilde \rho}

\begin{document}
\begin{center}
\thispagestyle{empty}

{\Large Conditional Limit Results for Type I Polar Distributions}

       \vskip 0.4 cm

         \centerline{\large Enkelejd Hashorva}

        \vskip 0.8 cm

        \centerline{\textsl{Department of Mathematical  Statistics and Actuarial Science}}
        \centerline{\textsl{University of Bern, Sidlerstrasse 5}}
        \centerline{\textsl{CH-3012 Bern, Switzerland}}
\today{}

\end{center}

%\subjclass[2000]{Primary 60F05; Secondary ?}

%%%%%%%%%%%%%%%%%%%%%%%%%%%%%%%%%%%%%%%%%%%%%%
{\bf Abstract:}  Let $(S_1,S_2)=(R \cos(\Theta), R \sin (\Theta))$
be a bivariate random vector with associated random radius $R$
which has distribution function $F$ being further independent of the
random angle $\Theta$. In this paper we
investigate the asymptotic behaviour of the conditional survivor
probability $\overline{\Psi}_{\rho,u}(y):=\pk{\rho S_1+ \sqrt{1-
\rho^2} S_2> y \lvert S_1> u}, \rho \in (-1,1),\inr$ when $u$ approaches
the upper endpoint of $F$. On the density function of $\Theta$ we require a certain local
asymptotic behaviour at 0, whereas for $F$ we require that it belongs to the Gumbel max-domain of attraction. The main result of
this contribution is an asymptotic expansion of
$\overline{\Psi}_{\rho,u}$, which is then utilised to construct two
estimators for the conditional distribution function $1- \overline{\Psi}_{\rho,u}$. Further, we allow $\Theta$ to depend on $u$.


%*********************************************************************************************************

\section{Motivation}
Let $(S_1,S_2)$ be a spherical bivariate random vector with
associated random radius $R>0$ (almost surely) with
distribution function $F$. The random vector $(X,Y)$ with stochastic representation
\BQNY
 (X, Y) &\equaldis &(S_1, \rho S_1+ \sqrt{1- \rho^2} S_2), \quad \rho\in (-1,1)
 \EQNY
is an elliptical random vector ($\equaldis$ stands for equality of
the distribution functions). If $F$ is in the Gumbel max-domain of attraction with positive scaling function $w$, i.e.,
\BQN \label{eq:gumbel:w}
\lim_{u \uparrow x_F} \frac{1 - F(u+x/w(u))}{1- F(u)}&=& \exp(-x),
\quad \forall x\inr,
\EQN
where  $x_F\in (0, \IF]$ is the upper endpoint of $F$,  then Theorem 4.1 in Berman (1983) implies the following
Gaussian approximation
\BQN\label{eq:ber}
\lim_{u \uparrow x_F}\pb{Z_{u,\rho} >  \rho u+ y \sqrt{u/w(u)} }& =&
\pb{ Z> y/\sqrt{1- \rho^2}}, \quad \forall y\inr,
\EQN
with $Z_{u,\rho}\equaldis Y  \lvert X > u$ and $Z$ a standard Gaussian random variable (mean 0 and variance 1).

Berman's result shows that the Gumbel max-domain of attraction assumption is
crucial for the derivation of \eqref{eq:ber}. Conditional
limit results for $F$ in the Weibull max-domain of attraction and
$(X,Y)$ a bivariate  elliptical random vector are obtained in Berman
(1992), Hashorva (2007b). The case $F$  is in the \FRE max-domain of
attraction is simpler to deal with, see Berman (1992).

As shown in Cambanis et al.\ (1981) we have the following stochastic
representation \BQN \label{eq;manuel} (S_1,S_2) &\equaldis &(R
\cos(\Theta), R \sin(\Theta)), \EQN with $R$ independent of the
random angle $\Theta$ which is uniformly distributed on $(-\pi,
\pi)$, i.e., $(\cos(\Theta))^2$  possesses the Beta distribution with parameters $1/2,1/2$.

When  $(\cos(\Theta))^2$ is Beta distributed, then the random vector $(S_1,S_2)$ is a generalised symmetrised Dirichlet
random vector.  Generalisation of \eqref{eq:ber} for such $(S_1,S_2)$ is presented in Hashorva (2008c) with
limit random variable $Z$ being Gamma distributed (see below Example 1).

Three natural questions arise:\\
a) What is the adequate approximation of the conditional survivor
function $\pk{Z_{u,\rho}>y}$ if $\Theta \in (-\pi, \pi)$ is some general random angle with unknown distribution function?\\
b) What can be said about the limit random variable $Z$?\\
c) Does $Z$ has a more general distribution if the random angle $\Theta=\Theta_u$ varies with $u$?

In this paper we show that if $\Theta_u$ possesses a positive
density function $h_u$ with a certain local asymptotic behaviour at
0, then we can answer both questions raised above. The
generalisation of \eqref{eq:ber} for bivariate polar random vectors
(see Definition \ref{eq:def:pol} below) satisfying
\eqref{eq:gumbel:w} is given in Section 3. Two applications of our
results are presented in Section 4. The first one concerns the
asymptotic behaviour of survivor function of bivariate polar random
vectors. In the second application we discuss the estimation of the
conditional distribution function $\pk{Z_{u,\rho}>y}$.
   Proofs and related results are relegated to Section 5.

\section{Preliminaries}
We shall explain  first the meaning of some notation, and then we introduce the class of bivariate polar random vectors.
A set of assumptions needed to derive the main results of this paper concludes this section.

If $X$ is a random variable with distribution function $H$ this will be alternatively denoted by $X\sim H$. When
$H$ possesses the density function $h$ we write $X \backsimeq h$.

In the following $\psi$ is a positive measurable function such that for all $z\in (0,\IF)$
\BQN \label{eq:psi}
\psi(z)&\le& K \max(z^{\lambda_1}, z^{\lambda_2}), \quad K> 0, \lambda_i\in (-1/2,\IF),\quad i=1,2,
\EQN
where $\E{\psi(W^2/2)}> 0$ with $W\sim \Phi$. Since $\E{\psi(W^2/2)}< \IF$ we can define
a distribution function $\Psi$ on $\R$ by
\BQN\label{eq:Psi}
\Psi(z):=\frac{\int_{-\IF }^{z } \exp(-s^2/2)\psi(s^2/2)\, ds}
{\int_{-\IF}^{\IF } \exp(-s^2/2)\psi(s^2/2)\, ds},
\quad \forall z\inr. \EQN
We denote by $\Psi_{\alpha, \beta}, \alpha,\beta>0$ the Gamma distribution with density function
$x^{\alpha-1} \exp(- \beta x) \beta^\alpha/\Gamma(\alpha), x\in (0,\IF) ,$ where $\Gamma(\cdot)$ is the Gamma function.\\
Next, we introduce the class of bivariate polar random vectors. Throughout the paper $R$ denotes a positive random radius with distribution function $F$ independent of the random angle $\Theta \in (-\pi,\pi)$, and $(S_1,S_2)$ is a bivariate random vector with representation \eqref{eq;manuel}. In the special case $\Theta$ is uniformly distributed on $(-\pi, \pi)$ for any two constants $a_1,a_2$ (see Lemma 6.1 in Berman (1983)) we have
\BQN\label{eq:Berm:a1}
 a_1S_1+ a_2 S_2 \equaldis \sqrt{a_1^2 + a_2^2} S_1 \equaldis \sqrt{a_1^2 + a_2^2} S_2,
 \EQN
hence linear combinations of spherical random vectors (i.e.\ the elliptical random vectors) are very tractable.

If  the random angle $\Theta$ is not  uniformly distributed on $(-\pi, \pi)$,
then \eqref{eq:Berm:a1} does not hold in general. In this paper we do not make specific distributional assumptions on $\Theta$. We assume however that the random angle $\Theta$ possesses a positive density function $h$ on $(-\pi, \pi)$.

\begin{defn} \label{eq:def:pol}
A bivariate random vector $(X,Y)$ is referred to as a bivariate
polar random vector with coefficients $a_i,b_i, i=1,2$ if it has the
stochastic representation
\BQN \label{eq:def:polar} (X,Y) \equaldis
(a_1S_1+a_2S_2, b_1 S_1 + b_2S_2), \quad (S_1,S_2)\equaldis  (R
\cos(\Theta), R \sin(\Theta)), \EQN where $R\sim F$  and $R>0$ (almost surely) being
independent of the random angle $\Theta \in (-\pi, \pi)$.
\end{defn}

Clearly, bivariate elliptical random vectors are included in the above class, which is defined in terms of
three components, a) the distribution of the associated random radius $R$, b) the distribution function of the random angle $\Theta$, and c) the deterministic coefficients $a_1,a_2,b_1,b_2$. In this paper we consider for simplicity the case
$$a_1=1,a_2=0,\text{ and } b_1=\rho, b_2= \sqrt{1- \rho^2}, \quad \rho \in (-1,1).$$
We refer to $\rho$ as the pseudo-correlation coefficient, and call $(X,Y)$ simply a bivariate polar random vector with
pseudo-correlation coefficient $\rho$. We have thus the stochastic representation
\BQN \label{eq:def:polar:2}
(X,Y) \equaldis (S_1, \rho S_1 + \sqrt{1- \rho^2}S_2), \quad (S_1,S_2)\equaldis  (R \cos(\Theta), R \sin(\Theta)), \quad
R\sim F,
\EQN
with $R>0$ independent of $\Theta$.\\
We note in passing that $S_1,S_2$ are in general dependent random
variables. If $S_1$ and $S_2$ are independent,  for instance if $R^2$ is chi-squared distributed with 2
degrees of freedom and $\Theta$ is uniformly distributed on
$(-\pi,\pi)$, then  $(X,Y)$ is a linear combination of
independent Gaussian random variables.

Next, we formulate three assumptions needed in this paper:

{\bf A1}. [Gumbel max-domain of attraction]\\
The distribution function $F$ with upper endpoint $x_F$ is in  the Gumbel max-domain of attraction satisfying \eqref{eq:gumbel:w}
with the scaling function $w$. Further,  suppose  that $F(0)=0$ and $x_F\in (0,\IF]$.

We formulate next an assumption for the second order approximation in \eqref{eq:gumbel:w}
initially suggested in Abdous et al.\ (2008).

{\bf A2}. [Second order approximation of $F$]\\
Let $F$ be a distribution function on $[0,\IF)$ satisfying Assumption A1. Suppose that there
exist positive functions $A,B$ such that
\BQN\label{main:A}
\Abs{ \frac{ 1- F(u+ x/w(u))}{1- F(u)}- \exp(-x)} \le A(u)B(x)
\EQN
holds for all $u< x_F$ large enough and any $x\in [0, \IF)$. Furthermore we assume $\lim_{u \uparrow x_F} A(u)=0,$
 and $B$ is locally bounded on finite intervals of $[0,\IF)$.

{\bf A3}. [Local approximation of $h_n,n\ge 1$ along $t_n$]\\
Let $h_n:(-\pi, \pi)\to [0,\IF), n\ge 1$ be a sequence of density functions such that $h_n(\theta)=h_n(- \theta), \forall \theta \in [0, \pi/2)$, and let $t_n,n\ge 1$ be
 positive constants tending to $\IF$ as $n\to \IF$.
Assume that for  any sequence of positive  measurable functions $\tau_n(s)= 1+ O(s/t_n),n\ge 1, s\ge 0$
for all large $n$ we have
\BQN \label{eq:has}
h_n\biggl(\tau_n(s)\sqrt{\frac{2z}{t_n}}\biggr)& =& h_n(1/\sqrt{t_n})
\psi_{\tau_n}(z\tau_n(s)),\quad \forall s,z\in [0,\IF),
\EQN
where $\psi_{\tau_n},n\ge 1$ are positive
measurable functions such that
\BQNY
 \psi_{\tau_n}(s) &\to &\psi(s), \quad n\to\IF
\EQNY
locally uniformly for $s,z \in [0,\IF)$ with $\psi_{\tau_n} $ satisfying $\eqref{eq:psi}$ for all large $n$ and all $s\in [0,\ve t_n)$ with $\ve$ a fixed positive constant.

Next we impose an assumption on the second order asymptotic behaviour of $h_n,n\ge 1$ at $0$.

{\bf A4}. [Second order approximation of $h_n,n\ge 1$ along $t_n$]\\
Suppose that Assumption A3 holds for some given sequence $t_n,n\ge 1,$ and further for any sequence of functions
$\tau_n(s)= 1+ O(s/t_n),n\ge 1, s\ge 0$ for all large $n$  we have
\BQN \label{eq:has:b}
\Abs{  \frac{1}{h_n(1/\sqrt{t_n})} h_n\biggl( \tau_n(s)\sqrt{\frac{2z}{t_n}}\biggr)-  \psi(z) }\le a(t_n) b_n(z),
\quad \forall s, z\in [0,\IF)
\EQN
where  $a,b_n,n\ge 1 $ are positive measurable functions such that
 $$\limit{n} a(t_n)= 0, \quad \limit{n} b_n(s)=b(s),$$
 and $b_n,n\ge 1$ satisfy \eqref{eq:psi}  for all $n$ large. % and all $s\in [0,\ve t_n), \ve >0$. ?

\section{Main Results}
In this section we consider a bivariate polar random vector $(X,Y)$ with pseudo-correlation $\rho$ and  representation \eqref{eq:def:polar:2}. We are interested in  the asymptotic behaviour of the conditional
distribution $Y \lvert X> u_n$ when $u_n$ tends ($n\to \IF$) to the upper
endpoint $x_F$ of $F$. Several authors have dealt with such
conditional probabilities and their statistical estimation,
see e.g., Gale (1980), Eddy and Gale (1981), Berman (1982, 1983, 1992), Heffernan and Tawn (2004), Abdous et al.\ (2005),
Heffernan and Resnick  (2007), Abdous et al.\ (2008) and Hashorva (2008b,c). Statistical modelling of conditional distributions is treated in the excellent monograph
Thomas and Reiss (2007).

The main assumption imposed on $F$ is that it satisfies  Assumption A1 with the scaling function $w$. Such polar random vectors are referred to alternatively as Type I polar random vectors.
The scaling function $w$ possesses two crucial asymptotic properties: a) uniformly on the compact sets of $ \R$
 \BQN\label{eq:wb}
\lim_{u \uparrow  x_F} \frac{w(u+ z/w(u))}{w(u)}&=&1, \EQN
and b)
\BQN\label{eq:uv}
 \lim_{u \uparrow x_F } u w(u) &= &\IF , \quad \lim_{u \uparrow x_F } w(u) (x_F- u)  = \IF \text{  if } x_F < \IF.
\EQN
Refer to Falk et al.\ (2004) or Resnick (2008) for  more details on the Gumbel max-domain of attraction.\\

We derive in the next theorem the asymptotic behaviour of $R\cos(\Theta_n)$, with $\Theta_n$ a random angle depending on $n$.

\BT\label{theo:0} Let $R$ be a positive random radius with
distribution function $F$ independent of the random angle $\Theta_n \backsimeq h_n, n\ge 1$. Let $u_n, n\ge 1$ be constants such that $u_n< x_F, n\ge 1$ and  $ \limit{n}u_n =x_F$ with $x_F\in (0,\IF]$ the upper endpoint of $F$. If $F$ satisfies Assumption A1, and
the density functions $h_n,n\ge 1$ satisfy Assumption A3 along $t_n:= u_n w(u_n), n\ge 1$ with $ \psi, \psi_{\tau_n},n\ge 1$,
then we have
\BQN \label{eq:theo:0}
\pk{R\cos(\Theta_n)> u_n}&=&  (1+o(1))  t_n^{-1/2} h_n(1/\sqrt{t_n}) [1- F(u_n)]
\int_{-\IF}^{\IF } \exp(-x^2/2)\psi(x^2/2)\, dx, \quad n\to \IF.
\EQN
If $\Theta_n=\Theta,\forall n\ge 1$ not depending  on $n$, then  $R \cos(\Theta)$
has distribution function in the Gumbel max-domain of attraction with the scaling function $w$. Furthermore, the
convergence in probability \BQN\label{eq:abs}
 q_n \abs{R\cos(\Theta)- u_n} \Bigl \lvert R\cos(\Theta)> u_n &\toprob & 0, \quad n\to \IF
\EQN
holds for any sequence $q_n,n\ge 1$ such that $\limit{n} w(u_n)/q_n =\IF$.
\ET
We note in passing that \eqref{eq:theo:0} is obtained in Theorem 12.3.1 of Berman (1992) assuming that
$(\cos(\Theta_n))^2$ is Beta distributed with positive parameters $a,b$. See also Tang (2006, 2008) for some important results on
tail asymptotics of products of independent random variables.

We state now the main result of this section.

\BT\label{theo:1} Let $(X_n,Y_n),n\ge 1$ be a bivariate polar random vector
with representation \eqref{eq:def:polar:2}, where
$\rho\in (-1,1), R\sim F$ and $\Theta_n\backsimeq h_n,n\ge 1$. Let $u_n,n\ge 1$ be a positive sequence such that $u_n< x_F,n\ge 1$ and $\limit{n} u_n = x_F$. Suppose that $F$ satisfies
Assumption A1 and $h_n,n\ge 1$ satisfy Assumption A3 along $t_n:= u_nw(u_n),n\ge 1$ with $ \psi, \psi_{\tau_n},n\ge 1$. If further
$\limsup_{n\to \IF} h_n((1+o(1))/\sqrt{t_n})/h_n(1/\sqrt{t_n}) < \IF$, then for any $x>0, y\inr$ we have
\BQN\label{eq:light:1}
 \limit{n}  \pb{ Y_n \le \rho u_n +y u_n/ \sqrt{t_n}, X_n \le u_n +x/w(u_n) \Bigl \lvert X_n> u_n }&=&\pb{Z \le y/\sqrt{1- \rho^2},W\le x},
\EQN
with $Z \sim \Psi$ being independent of $W\sim \Psi_{1,1},$ where  $\Psi$ is defined in \eqref{eq:Psi}.
\ET
Assumption A3 is somewhat cumbersome.  If we consider random angles $\Theta_n$ not depending on $n$ for all large $n$,
a tractable condition on the local asymptotic behaviour of the density of $\Theta_n$ is imposed below.

\BT \label{theo:5} Under the setup of \netheo{theo:1} if $\Theta_n=\Theta\backsimeq h, n\ge 1$ and instead of Assumption A3
we suppose that the density function $h$ of $\Theta$ is regularly varying at $0$ with index $2\delta\in (-1,\IF)$,
then for any sequence $u_n< x_F,n\ge 1$ such that $\lim_{n\to \IF} u_n=x_F$ we have
\BQN \label{eq:theo:0:5}
\pk{X_n > u_n}&=&  (1+o(1)) \frac{2 ^{\delta+1/2}}{\Gamma(\delta+1/2)}t_n^{-1/2} h(\sqrt{1/t_n})  [1- F(u_n)],
 \quad n\to \IF,
\EQN
and $X_1$ has distribution function in the max-domain of attraction of the Gumbel  distribution with the scaling function $w$. Furthermore  \eqref{eq:abs}
holds for any sequence $q_n,n\ge 1$ such that $\limit{n} w(u_n)/q_n =\IF$, and
for $x>0,y\inr$  given constants \eqref{eq:light:1} is satisfied with $Z^2 \sim \Psi_{\delta+1/2,1/2},$ and
$Z$ symmetric about 0 independent of $W\sim \Psi_{1,1}$.
\ET
We present  next an illustrating example.\\

{\bf Example 1.} [Kotz Type III Polar Random Vector]
Let $R\sim F$ be a random radius with tail asymptotic behaviour
\BQN\label{eq:kotz:Fk}
1- F(u) &=&  (1+o(1))K  u^{N}\exp(-r u^\delta), \quad K>0,\delta>0, N\inr, \quad u\to \infty.
 \EQN
If  $\Theta \backsimeq h $ is a random angle independent of $R$ we call $ (X,Y)$ with stochastic representation \eqref{eq:def:polar:2} a Kotz Type III polar random vector with pseudo-correlation $\rho\in (-1,1)$. If we set  $w(u):= r \delta u ^{\delta-1},u>0$, then
\BQNY
\limit{u} \frac{ \pk{R> u+ x/w(u)}}{\pk{R>u}} &=&\exp(-x), \quad \forall x\inr
\EQNY
implying that $F$ is in the Gumbel max-domain of attraction with the scaling function $w$. Suppose that $h(\theta)=h(-\theta),\forall \theta \in [0,\pi/2)$, and further
$$ h(\theta)=
c_{a,b}\abs{\sin(\theta)}^{2a-1} \abs{\cos(\theta)}^{2b-1}, \quad \theta \in (-\ve,\ve), \quad \ve\in (0,\pi),$$
where $a,b,c_{a,b}$ are positive constants. Note that when
$$\ve=\pi, \quad c_{a,b}=\frac{1}{2} \frac{\Gamma(a+b)}{\Gamma(a)\Gamma(b)}, $$
then $(X,Y)$ is a generalised symmetrised Dirichlet random vector (see Hashorva (2008c)). It follows that Assumption A3 is satisfied with
$$ h(1/\sqrt{t_n})= (1+o(1))c_{a,b} t_n^{1/2-a}, \quad \psi(s)= (2s)^{a- 1/2},\quad s>0, \quad t_n \to \IF$$
and $h$ is regularly varying at 0 with index $2a-1$. By \eqref{eq:theo:0}  for $u_n\to \IF$ we have
\BQNY
\pk{X > u_n}&=&   (1+o(1)) c_{a,b} K (2/(r\delta))^{a}\Gamma(a) u_n^{N- a\delta} \exp(-r u_n^\delta).
\EQNY
Next, for any $x>0, y\inr$ \netheo{theo:1} implies
\BQNY
\lim_{n\to \IF}\pb{  Y \le \rho u_n+  y u_n^{1-\delta/2},
X\le  u_n+ x  u_n^{1-\delta} \Bigl \lvert X> u_n }
&=&\pk{Z\le y \sqrt{ r \delta/(1- \rho^2)},W\le r \delta x},
\EQNY
with $Z$ symmetric about 0 independent of $W\sim \Psi_{1,1}$, and  $Z^2 \sim\Psi_{a,1/2}$.
Remark that if $a=1/2$, then  $Z$ is a standard Gaussian random vector. When also $b=1/2$,
then $(X,Y)$ is an elliptical random vector with pseudo-correlation $\rho$.

In the next theorem we show a second order correction for the conditional limit result obtained in \eqref{eq:light:1}
which is of some interest for statistical applications.

\BT\label{theo:3} Under the assumptions and the notation of
\netheo{theo:1}, if furthermore  Assumptions A2 and A4 are satisfied where $x_F=\IF$ and $\rho\in [0,1)$,
 then we have locally uniformly for any $z\inr$ (set $z_{n,\rho}:=\rho u_n+
zu_n\sqrt{1- \rho^2}/ \sqrt{t_n}$)
\BQN \label{eq:theo:3:R}
\pb{Y_n>z_{n,\rho}\Bigl \lvert X_n> u_n} &=&
1- \Psi(z) + \frac{1}{\sqrt{t_n}}\frac{\rho}{\sqrt{1- \rho^2} }\Psi'(z)+  O\Bigl( A(u_n) + a(t_n ) +\frac{1}{t_n}\Bigr) , \quad \ntoi,
\EQN
provided that
$$\max\Biggl( \int_0^\IF B(s)ds , \int_0^\IF B(s)\max(s^{\lambda_1},s^{\lambda_2})ds \Biggr) < \IF,$$
where $\lambda_i\in (-1/2, \IF),i=1,2$ are the constants related to Assumption A3.
\ET

\begin{remark}
a) Abdous et al.\ (2008) show several examples of distribution functions  $F$ satisfying Assumption A2.
The assumptions on $h$ can be easily checked for common distribution functions using Taylor expansion.

b) If we assume $h$ is regularly varying with index $2 \delta \in (-1,\IF)$ instead of the Assumption A3 and
modifying A4 accordingly, then \eqref{eq:theo:3:R} holds with $\Psi:=\Psi_{\delta+1/2,1/2}$,  provided that
$$\max\Biggl( \int_0^\IF B(s)ds , \int_0^\IF B(s)s^{\delta}ds \Biggr) < \IF.$$
\end{remark}

\section{Applications}
In this section we present two applications of our asymptotic
results: a) we obtain an  asymptotic expansion for the joint
survivor probability of polar random vectors, and b) we discuss briefly the
estimation of the conditional distributions of such vectors.

\subsection{Tail Asymptotics}
Let $(X,Y)$ be a bivariate polar random vector with pseudo-correlation coefficient $\rho\in (-1,1)$.
Assume that the distribution function $F$ of the random radius $R$ has an infinite upper endpoint.
In various situations quantification of the asymptotics of the joint survivor probability $\pk{X>x,Y>y}$ is
of interest when $x,y$ become large. Our asymptotic results in Section 3 imply an asymptotic expansion of this survivor probability, provided that $(X,Y)$ is of Type I. Explicitly, under the assumptions of \netheo{theo:5} we obtain for any $x>0, y\inr$ and $u$ large
(set $x_u:=u+ x/w(u), y_{u,\rho}:= \rho u + y\sqrt{u /w(u)}, u>0$)
\BQNY
\pk{X> x_u, Y>y_{u,\rho}} &=& (1+o(1))\exp(-x)[1- \Psi_{\delta+1/2, 1/2}(y)] \pk{ X> u}, \quad u\to \IF.
\EQNY
In our asymptotic result the sequence $y_{u,\rho}$ increases like $\rho u$ since by (13)
$$ y_{u,\rho}= (1 + \frac{y}{ \rho \sqrt{u w(u)}})\rho u=(1+o(1))\rho u, \quad u\to \IF.$$
It is of some interest to consider also constants $y_{u,\rho} =c u, u>0,$ with $c\in (\rho,1]$. In view of \netheo{theo:5} for any
$c\in (-\IF, \rho)$ we have
\BQNY \pk{X> x_u, Y>y_{u,c}}&=&(1+o(1))\exp(-x)\pk{ X> u}, \quad u\to \IF.
\EQNY When $c\in (\rho, 1]$ the joint survivor probability $\pk{X> x_u, Y>y_{u,c}}$
diminishes faster than $\pk{X>u}$,
i.e., \BQNY \limit{n}\frac{\pk{X> x_u, Y>y_{u,c}}}{\pk{ X> u}}&=&0.
\EQNY
If $(X,Y)$ is a bivariate elliptical random vector we may write (see Hashorva (2007c))
\BQN \label{eq:ext} \pk{X> u, Y> c u}
&=& (1+o(1))\frac{\alpha_{\rho,c} K_{\rho,c}}{2 \pi}\frac{1- F(\alpha_{\rho,c} u)}{u
w(\alpha_{\rho,c} u)}, \quad u\to \IF
\EQN
for any $c\in (\rho,1]$ with
\BQNY
\alpha_{\rho,c}:=\sqrt{(1 - 2c\rho + \rho^2)/(1- \rho^2)}\in (1,\IF) ,
 \quad  K_{\rho,c}&:=&\frac{(1- \rho^2)^{3/2} }{(1- c\rho )(c- \rho )}\in (0,\IF).
 \EQNY
In a forthcoming paper we extend \eqref{eq:ext} to the case of Type I bivariate polar random vectors.

\subsection{Estimation of Conditional Distributions}
Let $ (X_i, Y_i), i\le n,n\ge 1 $  be independent and identically distributed bivariate polar random
vectors with pseudo-correlation coefficient $\rho \in (-1,1)$ and random radius $R\sim F$. Define the conditional distribution
function
$$\Psi_{\rho,x}(y):= \pk{Y_1\le  y \lvert X_1> x}, \quad x,y\inr.$$
Suppose that $x_F=\IF$ and $F$ satisfies Assumption A1. As in the elliptical setup (Abdous et al.\ (2008))
also in the general case of polar random vectors estimation of the conditional distribution function
$\Psi_{\rho,x}$ can be motivated by our novel limit results, since under the assumptions of \netheo{theo:5}
 we have (set $t_u:=u w(u),u>0$)
 \BQN\label{gaussA}
\sup_{y\inr} \Abs{\Psi_{\rho,u}( u[ \rho + y \sqrt{1/t_u}] ) - \Psi_{\delta+1/2,1/2}( y/\sqrt{1- \rho^2})}
\to  0, \quad u\to \IF,
\EQN
with $2\delta$  the index of the regular variation of $h$ at 0. Under Assumptions A2 and A4 we obtain additionally the second order asymptotic expansion
 \BQN\label{gaussA}
\Psi_{\rho,u}(u[ \rho+1/ t_u + y\sqrt{1/t_u}]) =
\Psi_{\delta+1/2,1/2}( y/\sqrt{1- \rho^2}
)+O(A(u)+a(u)+\frac{1}{t_u}), \quad u\to \IF. \EQN
These approximations motivate the following estimators of $\Psi_{\rho,x}$ for $x$ large and $y$ positive, namely
$$\hat \Psi_{\rho,x,n}^{(1)}(y):= \Psi_{\delta+1/2,1/2}\Bigl(  \frac{ y- \hat \rho_n x  }{
\sqrt{(1- \hat \rho_n^2)x/\hat w_n(x)}}\Bigr), \quad n> 1,
$$
and
$$\hat \Psi_{\rho,x,n}^{(2)}(y):= \Psi_{\delta+1/2,1/2}\Bigl(  \frac{ y- \hat \rho_n (x+1/\hat w_n(x))  }{
\sqrt{(1- \hat \rho_n^2)x/\hat w_n(x)}}\Bigr), \quad n> 1,
$$
where  $\hat \rho_n$ is an estimator of $\rho$, and $\hat w_n(\cdot)$ is an estimator of $w(\cdot)$.

An estimator of $\hat \rho_n$ can be constructed considering the relation between $\rho$ and
the expectation $\E{Y}$, provided that the latter exists.
Estimation of $\delta$ and $w$ are difficult tasks. If the scaling function $w$ (related to
the Gumbel max-domain of attraction of $F$) is simple, say
$w(u)= c \gamma u^{\gamma-1}, c,\gamma>0, u>0,$ then an estimator $\hat w_n$ is constructed by
estimating separately $c$ and $\gamma$ from  $X_1 \ldot X_n$ (recall $X_1$ has distribution function in the Gumbel max-domain
of attraction with the scaling function $w$). See Abdous et al.\ (2008), Hashorva (2008a) for more details.

In practical situations also the constant $\delta$ might be unknown and  therefore  has  to be estimated.
One possibility of estimating $\delta$ is to utilise  \eqref{eq:theo:0:5}.

We note that for elliptical random vectors $\delta=2$ and both estimators $\hat \Psi_{\rho,x,n}^{(1)}$ and $\hat \Psi_{\rho,x,n}^{(2)}$
are suggested in Abdous et al.\ (2008). Since we estimate both $c$ and $\gamma$
from $X_1 \ldot X_n$  and not from the observations of the random radius $R$, our estimators above differ from those in the aforementioned paper.

\section{Related Results and Proofs}
Set in the following
\BQN
\aroxy:= \sqrt{1+ ((y/x)- \rho)^2/(1- \rho^2)}, \quad \broxy:= \aroxy x/y, \quad x,y\inr,y\not=0,\quad  \rho \in (-1,1).
\EQN
For $1 \le  a<  b\le \IF, x>0$ constants, and $h,F$ two positive measurable functions we define
\BQN
\label{iax} J(a,b, x,h):= \int_a^b [1- F(x t)] h(t) \frac{1}{t \sqrt{t^2- 1}}\, dt.
\EQN
If $b=\IF$ write simply  $J(a,x,h)$  suppressing the second argument. Write $\wth (\cdot)$ and $\btrho (\cdot) $ instead of $h(\arccos(1/\cdot)$ and $h(\arcsin(1/\cdot)-\arcsin(\rho))$, respectively.

Next, we shall prove two lemmas, and then proceed with the proof of
the main results. The first lemma is formulated for $F$ with
infinite upper endpoint. It generalises Lemma 5 in Hashorva (2008b)
for bivariate elliptical random vectors. If $F$ has a finite upper
endpoint, say $x_F=1$, then a similar result holds. Statement $b)$
and $c)$ should be reformulated requiring additionally that $x^2+2
\rho xy+ y^2< 1- \rho^2$ with $\abs{x},\abs{y}\in [0,1]$.

%%%%%%%%%%%%%%%%%%%%%%%%%%%%%%%%%%%%%%%%%%%%%%%%%%%%%%%%%%

\BL \label{lem:x1}
Let the random radius $R\sim F$ be independent of the random angle $\Theta \in (-\pi, \pi)$ and define a bivariate
polar random vector $(X,Y)$ with pseudo-correlation $\rho \in (-1,1)$ via \eqref{eq:def:polar:2}. If the upper endpoint $x_F$ of $F$ is infinite and $\Theta$ possesses a density function $h$ such that
$h(\theta)= h(-\theta), \theta \in [0,\pi/2)$, then we have:\\
a) For any $x> 0$
\BQN\label{eq:X}
\pk{X > x}&=& 2 J(1,x, \wth).
\EQN
b) For any $x>0,y\in (0,x]$ such that $y/x > \rho$
\BQN \label{eq:lem:A:A}
 \pk{X > x, Y> y}&=& J( \aroxy,x,\wth )+ J(\broxy,y,\btrho).
\EQN
c) For any $x>0$ and $y/x \in (0,\rho), \rho>0$
\BQN\label{eq:lem:A:B}
 \pk{X > x, Y> y}&=& 2 J(1,x,\wth) - J( \aroxy,x,\wth)+ J(\broxy ,y,\btrho).
\EQN
\EL
\begin{proof}
Since the associated random radius $R$ is almost surely positive being further independent of $\Theta$ and
$h(-\theta)= h(\theta), \cos(-\theta)= \cos(\theta), \forall \theta \in [0,\pi/2)$ for any $x> 0$ we obtain
\BQNY
\pk{X > x}&=& 2 \int_{0} ^{\pi/2}\pk{ R > x / \cos (\theta)} h(\theta)\, d \theta\\
&=& 2\int_1^\IF [1 - F(xs)] \frac{h(s)}{s \sqrt{s^2-1}}\, ds= 2 J(1,x, \wth).
\EQNY
We prove next the second statement. By the assumptions  $(X,Y) \equaldis  (R\cos(\Theta), R \sin(\Theta + \arcsin(\rho)))$.
Consequently for $x>0,y>0$  two positive constants
\BQNY
\pk{S_1> x, \rho S_1+ \sqrt{1- \rho^2} S_2> y}&=& \pk{R
\cos (\Theta) > x, R \sin (\Theta+ \arcsin(\rho)) > y}.
\EQNY
Since $\sin(\arcsin (\rho)+ \theta)/\cos(\theta)$  is strictly increasing
in $\theta \in [- \arcsin(\rho), \pi/2]$ with inverse $\arctan((\cdot- \rho)/\sqrt{1- \rho^2)}$ (see Kl\"uppelberg et al.\ (2007))
we have
\BQNY
%\lefteqn{\pk{R \cos (\Theta) > x, R \sin (\Theta+ \arcsin(\rho)) > y}}\\
\pk{X  > x, Y> y}&=&
\int_{\arctan((y/x- \rho)/\sqrt{1- \rho^2})}^{\pi/2} \pk{ R> x/ \cos(\theta)} \, d Q(\theta)\\
&&+ \int_{- \arcsin(\rho)}^{\arctan((y/x- \rho)/\sqrt{1- \rho^2})}
\pk{ R> y/ \sin(\theta+ \arcsin (\rho))} \, d Q(\theta),
\EQNY
with $Q$ the distribution function of $\Theta$. Transforming the variables we obtain for  $y/x > \rho$
\BQNY
\pk{X  > x, Y> y}&=& J(\aroxy,x,\wth)+ J(\broxy ,y,\btrho),
\EQNY
and if $y/x \le  \rho$ with $x,y$ positive
\BQNY
%\pk{R \cos (\Theta) > x, R \sin (\Theta+ \arcsin(\rho))> y}
\pk{X  > x, Y> y}&=&
2J(1,x,\wth) - J( \aroxy,x,\wth)+  J(\broxy ,y,\btrho),
\EQNY
hence the proof is complete.\end{proof}

\BL \label{lem:x2} Let $F$ be  a distribution function with upper endpoint $x_F\in (0,\IF]$ satisfying further
\eqref{eq:gumbel:w} with the scaling function $w$ and let $ 1 \le a_n\le  b_n, \gamma_n>1, u_n\in (0,x_F), t_n:=u_n w(\gamma_n u_n), n\ge 1$
be positive constants such that
\BQN
b_n u_n< x_F, n\ge 1, \quad \limit{n} \gamma_n= \gamma \in [1, \IF), \quad \limit{n} \gamma_n u_n =
\limit{n} b_n u_n= x_F,
\EQN
 and further
\BQN
\limit{n} t_n (a_n- \gamma_n)= \xi\in [0,\IF), \quad \limit{n} t_n (b_n- \gamma_n)&=&
\eta\in [\xi, \IF].
\EQN
Let $h,r,\psi_n,n\ge 1$  be positive measurable functions. Assume that for some $\ve>0$
\BQN \label{eq:has:b}
h(\gamma_n + s/t_n)& =&  r(\gamma_n,t_n) \psi_n(s), \quad \forall s \in [0, \ve t_n]
\EQN
and
\BQN \label{eq:has:2}
\psi_n(s) & \to & \psi(s) \in [0,\IF), \quad n\to \IF
\EQN
locally uniformly with $\psi_n$ satisfying \eqref{eq:psi} for all $n\ge 1, s\in [0, \ve t_n]$ with
$\lambda_i,i=1,2\in (c,\IF)$. Suppose further $\int_{a_n}^\IF h(s) (s \sqrt{s^2-1})^{-1}\, ds  < K< \IF , \forall n>1$. \\
a) If $\gamma \in (1, \IF)$ and $c=-1$
\BQN
 J(a_n, b_n, u_n, h)&=&
  (1+o(1)) \frac{ r(\gamma_n, t_n)}{ \gamma \sqrt{\gamma^2- 1}}
  \frac{1- F(\gamma_n u_n)}{t_n} \int_{\xi }^{\eta } \exp(-x)\psi(x)\, dx,   \quad n\to \IF.
\EQN
b) When  $\gamma=1$ and $\limit{n}  t_n (\gamma_n- 1)= \tau\in [0,\IF)$, then
\BQN
 J(a_n, b_n, u_n, h)&=&
  (1+o(1))r(\gamma_n,t_n)\frac{1- F(\gamma_n)}{\sqrt{t_n}  } \int_{\xi }^{\eta } \exp(-x)\frac{1}{\sqrt{2 x+ 2 \tau}}\psi(x)\, dx
  \quad n\to \IF,
\EQN
provided that $c=-1/2$ if $\xi =\tau=0$ and $c=-1$ otherwise.
\EL

\begin{proof} Set in the following for $n\ge 1$
$$ u_n^*:= \gamma_n u_n, \quad t_n:= u_n w(u_n^*), \quad l_n(x)= \gamma_n+ x/t_n, \quad
\psi_n^*(x):=\frac{\gamma\sqrt{\gamma^2- 1}}{l_n(x)
 \sqrt{l_n^2(x)- 1}} \psi_n(x), \quad x\ge 0$$
 and $ \xi_n:= t_n(a_n -\gamma_n), \eta_n:= t_n(b_n- \gamma_n).$ Since $\limit{n} u_n^*= x_F$, then \eqref{eq:uv}
 implies $ \limit{n} t_n = \IF, w(u_n^*)(x_F- u_n^*)= \IF.$ If  $c_1,c_2$ are two arbitrary constants such that $c_2>c_1>1$ for all $n$ large we have
\BQNY
\int_{\gamma+ c_2}^\IF [1- F(u_ns)] h(s)\frac{1}{s\sqrt{s^2-1}}\, ds
 &\le & [1- F(u_n(\gamma+ c_2))]\int_{\gamma+ c_2}^\IF  h(s)\frac{1}{s\sqrt{s^2-1}}\, ds
\EQNY
and
\BQNY
 J(a_n, b_n, u_n, h)&\ge &\int_{a_n}^{b_n} [1- F(u_ns)] h(s)\frac{1}{s\sqrt{s^2-1}}\, ds \\
 &\ge & \int_{a_n}^{\gamma+ c_1} [1- F(u_ns)] h(s)\frac{1}{s\sqrt{s^2-1}}\, ds \\
 &\ge & [1- F(u_n(\gamma+ c_1))] \int_{a_n}^{\gamma+ c_1} h(s)\frac{1}{s\sqrt{s^2-1}}\, ds.
 \EQNY
Assume that $x_F=\IF$. Since $1- F$ is rapidly varying (see e.g., Resnick (2008)) i.e.,
\BQNY
 \limit{n} \frac{1- F(u_nx)}{1- F(u_n)}&=&0 ,\quad \forall x>1
 \EQNY
for any $\ve^* >0$ we obtain
\BQNY
  J(a_n, b_n, u_n, h)&= & (1+ o(1))\int_{a_n}^{\gamma+ \ve^*} [1- F(u_ns)] h(s)\frac{1}{s\sqrt{s^2-1}}\, ds.
 \EQNY
If  $\gamma\in (1,\IF)$, then
$$ \psi_n^* (s) \to \psi(s),  \quad n\to \IF$$
locally uniformly for $s\ge 0$ and $\psi_n^*$ satisfying  \eqref{eq:psi} for all $s\in [0, \ve t_n), \ve>0$.
As in the proof of Lemma 7  of Hashorva (2007a) for any $\ve >0$ as $n\to \IF$ we obtain
\BQNY
\lefteqn{ J(a_n, b_n, u_n, h)}\\
&=& (1+o(1)) \int_{\gamma_n + \xi_n/t_n }^{\min(\gamma_n + \eta_n /t_n, \gamma+ \ve) } [1- F(u_n s)]
 h(t) \frac{1}{s \sqrt{s^2- 1}}\, ds\\
&=& \frac{(1+o(1))}{t_n} \int_{\xi_n }^{\min(\eta_n, t_n(\gamma- \gamma_n+ \ve) ) } [1- F(u_n^*+ x/w(u_n^*) )]
  h(l_n(x)) \frac{1}{l_n(x) \sqrt{l_n^2(x)- 1}}\, dx\\
&=& \frac{(1+o(1))}{\gamma \sqrt{\gamma^2- 1}} \frac{r(\gamma_n,t_n)}{t_n}
\int_{\xi_n }^{\min(\eta_n, t_n(\gamma- \gamma_n+ \ve) ) }[1- F(u_n^*+ x/w(u_n^*) )]
 \psi_n^* (x) % \frac{\sqrt{\gamma^2_n- 1}}{(1 + x/(t_n \gamma_n)) \sqrt{(1+ x/(t_n \gamma_n))^2- 1/\gamma_n^2}}
 \, dx\\
&=& (1+o(1)) %\frac{1}{\gamma \sqrt{\gamma^2- 1}}
\frac{1}{\gamma \sqrt{\gamma^2- 1}} \frac{r(\gamma_n,t_n)}{t_n }  [1- F(u_n^*)] \int_{\xi }^{\eta } \exp(-x)  \psi(x) \, dx.
\EQNY
Next, if $\gamma=1$ redefine
$$ \psi_n^*(s):= \frac{1 }{\sqrt{t_n}l_n(s) \sqrt{l_n^2(s)- 1}}\psi_n(s), \quad n\ge 1, s\ge 0.$$
We have
$$ \psi_n^*(s)\to  \frac{\psi(s) }{\sqrt{ 2\tau + 2 s}}=:  \psi^*(s)$$
locally uniformly for $s\ge 0$. Hence as in the proof above for $\ve>0$ and  $n\to \IF$ we obtain
\BQNY
J(a_n, b_n, u_n, h)&=&  (1+o(1))\frac{1}{\sqrt{t_n}}
\int_{\xi_n }^{\min(\eta_n, \ve t_n ) }
[1- F(u_n^*+ x/w(u_n^*) )]\psi^*_n(x)\, dx\\
&=& (1+o(1)) r(\gamma_n,t_n)  \frac{1- F(u_n)}{\sqrt{t_n}  }
\int_{\xi }^{\eta } \exp(-x)\frac{1}{\sqrt{2 x+ 2 \tau}}\psi(x)\, dx.
\EQNY
Similarly,  the asymptotic results  follow when $x_F \in (0,\IF)$,  hence the proof is complete. \end{proof}

\prooftheo{theo:0}
\COM{A2 implies (see Falk et al.\ (2004))
$ \limit{n} u_n w(u_n)= \IF$. If the upper endpoint $x_F$ of $F$ is finite, then we have further
\BQN\label{tn2}
\limit{n} w(u_n)(x_F- u_n)&=& \IF.
\EQN
The latter asymptotics is useful for the case $x_F$ is finite.}
We consider for simplicity only the case $x_F= \IF$. For all $n$ large \eqref{eq:X} implies
\BQNY
\pk{R \cos(\Theta_n) > u_n} &=& 2 \int_1^\IF [1- F(u_n s)] h_n( \arccos(1/x)) \frac{1}{x} \frac{1}{ \sqrt{x^2- 1}}\, dx.
\EQNY
We have (set $t_n:= u_n w(u_n),n\ge 1$)
\BQNY
\arccos(1/(1+ s/ t_n))
=\frac{\sqrt{2 s}}{\sqrt{t_n}} (1+ O(s/t_n))=: \sqrt{2 s/t_n}\tau_n(s), \quad n\to \IF
\EQNY locally uniformly for $s\ge 0$. Hence the Assumption A3 on $h_n$ implies
\BQNY h_n(\arccos(1/(1+ s/ t_n)))& = & h_n(  \tau_n(s)\sqrt{ 2 s/ t_n} )=h_n(1/\sqrt{t_n}) \psi_{\tau_n} (s\tau_n(s)), \quad s\ge 0.
 \EQNY
Applying \nelem{lem:x2} with $\tau-1=\gamma=\gamma_n=a_n=1,n\ge 1$ and $b_n=\IF, n\ge 1$ we obtain
\BQNY
\pk{R \cos(\Theta_n)> u_n}&=& (1+o(1))h_n(1/\sqrt{t_n}) \frac{1- F(u_n)}{\sqrt{t_n}  } \int_{-\IF}^{\IF }
\exp(-s^2/2)\psi(s^2/2)\, ds. \EQNY
If $h_n=h,n\ge 1$, then by the Assumption A3 we have $\limit{n}
h(1/\sqrt{t_n})/h(y_n/\sqrt{t_n})=1$ for any sequence $y_n,n\ge 1$ such that $\limit{n} y_n=1$. Consequently, the self-neglecting
property of $w$ in \eqref{eq:wb} implies
\BQNY \limit{n}
\frac{\pk{R \cos(\Theta) > u_n+ x/w(u_n)}}{\pk{R \cos(\Theta) > u_n}}&=& \limit{n} \frac{1-
F(u_n+ x/w(u_n)) }{1- F(u_n)} = \exp(-x), \quad \forall  x\inr.
\EQNY
Hence for any $z>0$
\BQNY \pk{
q_n \abs{R \cos(\Theta)- u_n}> z\lvert R \cos(\Theta)> u_n} &=&
\frac{\pk{ R \cos(\Theta)>  u_n + z/q_n}}{\pk{R \cos(\Theta)> u_n}}\\
&=&
\frac{\pk{ R \cos(\Theta)>  u_n + (z/w(u_n))(w(u_n)/ q_n)}}{\pk{R \cos(\Theta)>  u_n}}\\
&\to& 0,\quad n\to \IF,
 \EQNY
 thus the result follows. \QED

\prooftheo{theo:1}
 Set for $n\ge 1$ and $z\inr$
$$ v_n= z\sqrt{u_n/w(u_n)}, \quad \chi_n:= v_n/u_n, \quad  \arhon:=\sqrt{1+ (v_n/u_n)^2}, \quad t_n:= u_n w(u_n),\quad n\ge 1$$
and write in the sequel $\wthN (\cdot)$ and $\btrhoN (\cdot) $ instead of $h_n(\arccos(1/\cdot)$ and $h_n(\arcsin(1/\cdot)-\arcsin(\rho))$, respectively.

Since $\limit{n} t_n=\IF$ by the assumptions on $h$ making use of \eqref{eq:uv} and \eqref{eq:abs}
we retrieve the convergence in probability
$$ \sqrt{w(u_n)/u_n} (X_n- u_n) \Bigl \lvert X_n> u_n \toprob 0, \quad n\to \IF.$$
Consequently, it suffices to show the proof for $\rho=0$. Next, we prove the convergence in distribution
\BQNY \sqrt{w(u_n)/u_n} Y_n^* &\todis&  Z\sim \Psi, \quad n\to \IF,
\EQNY
with $Y_n^*\equaldis Y_n  \lvert X_n > u_n$ and $\Psi$ defined in \eqref{eq:Psi}.
Since  $\chi_n= v_n/u_n > \rho=0$ holds for all large $n$, we have in view of \nelem{lem:x1}
for all large $n$
\BQNY
\pk{X_n> u_n, Y_n>  v_n}
&=& J(\arhon,u_n,\wth )+ J(\chi^{-1}_n\arhon, v_n,\btrho ),
%.%\\
% &=&J(1+ (1+o(1))z_n^2/(2 u_n w(u_n)),u_n,\wth )+ J(((1+o(1))z_n^2/(2 u_n w(u_n)))/a_n , v_n, \btrho ).
 \EQNY
where $\arhon = 1+ (1+o(1))z^2/(2 t_n), n\to \IF.$ As in the proof of \netheo{theo:0} we obtain for the first term
\BQNY J(\arhon,u_n,\wth) &=& \frac{ (1+o(1))h_n(1/\sqrt{t_n})}{\sqrt{t_n}}
 [1- F(u_n)] \int_{z}^{\IF } \exp(-x^2/2)\psi(x^2/2)\, dx.
 \EQNY
Further, for any $s\ge 0$ (set $l_n(s):=\chi^{-1}_n+ s/(v_n w(u_n))$ we have
\BQNY
\frac{1}{l_n(s)}&=& \frac{\chi_n}{1 + s/t_n } =\sqrt{z^2/t_n}\frac{1}{1+ s/t_n}, \quad n\to \IF.
\EQNY
Consequently, the assumption on $h$ implies for all $s \ge 0$
\BQNY
\btrhoN ( l_n(s)) &=& h_n(\frac{z }{ \sqrt{t_n}}\tau_n(s)))= h_n(1/\sqrt{t_n})\psi_{\tau_n}\bigl(\tau_n(s)z^2/2\bigr), \EQNY
where $\tau_n(s):= 1+ O(s/t_n),s\ge 0,n\ge 1.$ Hence
\BQNY
\frac{1}{ l_n(s) \sqrt{ ( l_n(s))^2- 1}} \btrhoN ( l_n(s)) &=&
h_n(1/\sqrt{t_n}) \chi_n^{3/2} \psi_{\tau_n}\bigl(\tau_n(s)z^2/2\bigr), \quad \ntoi.
\EQNY
As in the proof of \nelem{lem:x2} we have thus \BQNY J(\chi_n^{-1} \arhon,
v_n,\btrho )
&=& \int_{\chi_n^{-1}\arhon }^\IF [1- F(v_n t)]\btrhoN(t)\frac{1}{t\sqrt{t^2- 1}}\, dt\\
&= & h_n(1/\sqrt{t_n})\frac{1- F(u_n)}{v_n w(u_n)} \chi_n^{3/2}
\int_{t_n[ \arhon- 1]}^\IF
\frac{1- F(v_n l_n(s))}{1- F(u_n)} \psi_{\tau_n}\bigl(\tau_n(s)z^2/2\Bigr)  \, ds\\
&= & (1+ o(1))h_n(1/\sqrt{t_n})\frac{1- F(u_n)}{v_n w(u_n)} \chi_n^{3/2} \psi(z^2/2)\int_{ z^2/2(1+o(1))}^\IF
\frac{1- F(u_n+ s/ w(u_n))}{1- F(u_n)}   \, ds\\
&=& o(J(1,u_n, \wthN)), \quad n\to \IF
\EQNY
implying
\BQNY
\limit{n} \pk{ Y_n^*> z \sqrt{u_n/w(u_n)} }&=& \limit{n}
\frac{ \pk{Y_n> z \sqrt{u_n/w(u_n)}, X_n> u_n}}{\pk{X> u_n}}\\
&=& \frac{\int_{z}^{\IF } \exp(-x^2/2)\psi(x^2/2)\, dx}{\int_{-\IF}^{\IF } \exp(-x^2/2)\psi(x^2/2)\, dx}=1- \Psi(z).
\EQNY
Thus the proof is complete. \QED

\prooftheo{theo:5} By the assumption on $h$ we have $h(s)= s^{2 \delta} L(s)$ for all $s>0$ in a neighbourhood of $0$
with  $L(s)$ a positive slowly varying function such that $\lim_{t \to 0} L(ts)/L(t)=1, \forall s>0$.
Furthermore, by Proposition B.1.10 in de Haan and Ferreira (2006) we have for any $\ve>0,\xi>0$
\BQN\label{eq:dehaan}
 \Abs{ \frac{h(t s )}{h(t)}- s^{2\delta}} \le \ve \max(s^{2 \delta- \xi},s^{2 \delta+\xi})
 \EQN
holds for any $s \in (0,t_0(\ve,\xi)/t), t\in (0,1)$ with $t_0(\ve,\xi)$ some positive constant. Since for positive constants $t_n,n\ge 1$ such that $\lim_{n \to \IF} t_n=\IF$
$$\frac{h(\sqrt{2s /t_n})}{h(\sqrt{1/t_n})}= (2s)^{\delta}=:\psi(s), \quad \forall s>0 $$
the result follows along the lines of the proof of \netheo{theo:1} utilising further \eqref{eq:dehaan}. \QED

\prooftheo{theo:3} Set for $n\ge 1, z\inr $ and $\abs{\rho}< 1$
$$ v_n:= \rho u_n+ z \sqrt{1- \rho^2}\sqrt{u_n/ w(u_n)}, \quad \chi_n:= v_n/u_n, \quad t_n:= u_n w(u_n),
\quad \arhon:= \sqrt{1 + (1/\chi_n^2- \rho)/(1- \rho^2)}.$$
In view of \nelem{lem:x1} for all large $n$ we have $ \pk{X_n> u_n}=  2J(1,u_n, \wthN)$
and further for $\rho\ge 0,z\ge 0$
$$ \pk{Y_n> v_n \lvert X_n> u_n}=  \frac{1}{ 2J(1,u_n, \wthN)} \biggl[J(\arhon,u_n, \wthN) +
J(\chi_n^{-1} \arhon, v_n, \btrhoN)\biggr].$$
In order to show the proof we need to approximate $J(1+ \frac{z^2}{2 t_n},u_n, \wthN)$ and
$J(\chi^{-1}_n\alpha_{n}, v_n, \btrhoN)$. We have
$$\arhon=1+ \frac{z^2}{2 t_n}+ O(1/t_n^2), \quad n\ge 1.$$
As in the proof of \netheo{theo:1} we obtain
\BQNY
 J(\arhon,u_n, \wth)
&=& (1+o(1)) h_n(1/\sqrt{t_n}) \frac{1- F(u_n)}{\sqrt{t_n}} [1- \Psi(z)]\int_{-\IF}^\IF
\exp(-x^2/2) \psi(x^2/2) \, dx,\quad n\to \IF.
\EQNY
Assumptions A3 and A4 imply $ \psi(x) \le \max(x^{\lambda_1}, x^{\lambda_2})$ and $b(x) \le \max(x^{\lambda_1^*}, x^{\lambda_2^*})$
 for some $\lambda_i, \lambda_i^*\in (-1/2,\IF)$.  Consequently, the assumptions on $B$ and $b$ yield
$$ \int_0^\IF \exp(-x) \psi(x) \sqrt{x}\, dx <  \IF, \quad \int_0^\IF B(x)\psi(x)\frac{1}{\sqrt{2x}}\, dx < \IF, \quad
\int_0^\IF \exp(-x) b(x) \frac{1}{\sqrt{2x}}\, dx <  \IF.$$
Define for all $x\ge 0$ and $n\ge 1$
$$ \gnt:=\frac{1}{(1+ x/t_n) \sqrt{t_n(1+ x/t_n)^2- t_n}}.$$
For all $x> 0, n\ge 1$ we have $\gnt \le \frac{1}{\sqrt{2x}}$, and further
\BQN
\Abs{\frac{1}{\sqrt{2x}}- \gnt}&\le &\frac{1}{\sqrt{2x}}\Abs{1-  \frac{1}{1+ x/t_n}}+ \frac{1}{1+ x/t_n}\Abs{ \frac{1}{\sqrt{2x}} - \frac{1}{ \sqrt{2x+ x^2/t_n}}}\le   5/4 \sqrt{x}/(t_n \sqrt{2}).
\EQN
Consequently, for $\zeta\ge 0$ and $n$ large we obtain
\BQNY
\lefteqn{\Abs{ \frac{\sqrt{t_n}}{ h_n(1/\sqrt{t_n})[1- F(u_n)]} J(1+ \zeta/ t_n,u_n, \wthN)- \int_{\zeta}^\IF
\exp(-x) \psi(x) \frac{1}{\sqrt{2 x}}\, dx }}\\
&= & \Abs{\int_{1+ \zeta/t_n}^\IF \frac{1- F(u_n t)}{1- F(u_n)} \frac{\sqrt{t_n}}{h_n(1/\sqrt{t_n})}
\wthN(t)\frac{1}{t \sqrt{t^2- 1}}\, dt - \int_{\zeta}^\IF
\exp(-x) \psi(x) \frac{1}{\sqrt{2 x}}\, dx }\\
&= & \Abs{\int_{\zeta}^\IF \frac{1- F(u_n+  x/w(u_n))}{1- F(u_n)} \frac{\wthN(1+ x/t_n)}{h_n(1/\sqrt{t_n})}
\gnt\, dx - \int_{\zeta}^\IF
\exp(-x) \psi(x) \frac{1}{\sqrt{2 x}}\, dx }\\
&\le &  \int_{\zeta}^\IF  \Abs{\frac{1- F(u_n+  x/w(u_n))}{1- F(u_n)} - \exp(-x)}
\frac{\wthN(1+ x/t_n)}{h_n(1/\sqrt{t_n})}\gnt\, dx \\
&&+\int_{\zeta}^\IF  \exp(-x) \Abs{  \frac{\wth(1+ x/t_n)}{h_n(1/\sqrt{t_n})}
- \psi(x) \frac{1}{\sqrt{2 x}}}\, dx \\
&\le &  A(u_n) \int_{\zeta}^\IF  B(x)\psi_n(x) \frac{1}{\sqrt{2x}} \, dx + \int_{\zeta}^\IF  \exp(-x)
\gnt\Abs{  \frac{\wth(1+ x/t_n)}{h_n(1/\sqrt{t_n})}- \psi(x)}\, dx \\
&&+\int_{\zeta}^\IF  \exp(-x) \psi(x) \Abs{ \gnt- \frac{1}{\sqrt{2 x}}}\, dx \\
&\le &  A(u_n)\int_{\zeta}^\IF  B(x)\psi(x)\frac{1}{\sqrt{2x}} \, dx
+ a(t_n) \int_{\zeta}^\IF  \exp(-x)\frac{1}{ \sqrt{2 x}} b_n(x)\, dx
+\frac{5\sqrt{2}}{8t_n} \int_{\zeta}^\IF  \exp(-x) \psi(x) \sqrt{x} \, dx \\
&= &  O\Bigl( A(u_n) + a(t_n)+ 1/t_n\Bigr)=:R_n(u_n).
\EQNY
Since
\BQNY
\Abs{\int_{z^2/2+ O(1/t_n)}^\IF
\exp(-x) \psi(x) \frac{1}{\sqrt{2 x}}\, dx - \int_{z^2/2}^\IF
\exp(-x) \psi(x) \frac{1}{\sqrt{2 x}}\, dx }
&= & O(1/t_n), \quad n\to\IF
\EQNY
we have
\BQNY
\Abs{ \frac{\sqrt{t_n}}{ h_n(1/\sqrt{t_n})[1- F(u_n)]} J(\arhon,u_n, \wthN)- \int_{z^2/2}^\IF
\exp(-x) \psi(x) \frac{1}{\sqrt{2 x}}\, dx }&= &  R_n(u_n).
\EQNY
Assume for simplicity in the following that $\rho >0$ and $z\ge 0$.
The other case can be established as in the proof of \netheo{theo:1}. We obtain the first order asymptotic expansion
(set $l_n(s):= \chi_n^{-1}+ s/(v_n w(u_n)), s>0,n>1$)
\BQNY
J(\chi^{-1}_n \arhon, v_n,\btrhoN ) &=& \int_{\arhon u_n/v_n}^\IF [1- F(v_n t)]\btrhoN(t)\frac{1}{t\sqrt{t^2- 1}}\, dt\\
&=&\frac{1- F(u_n)}{v_n w(u_n)}\int_{ v_n w(u_n)[\chi_n^{-1}\arhon- \chi_n^{-1}] }^\IF
\frac{1- F(u_n+ s/w(u_n))}{1- F(u_n)}
\btrhoN(l_n(s)) \frac{1}{l_n(s)\sqrt{l_n^2(s)- 1}}\, ds\\
&=&\frac{\rho^2}{ \sqrt{1- \rho^2}}h_n(1/\sqrt{t_n})\frac{1- F(u_n)}{v_n w(u_n)}\int_{ t_n[\arhon- 1]}^\IF
\frac{1- F(u_n+ s/w(u_n))}{1- F(u_n)} \psi_{\tau_n}(z^2/2(1+o(1)))\, ds\\
&=  & \frac{\rho^2}{ \sqrt{1- \rho^2}}h_n(1/\sqrt{t_n})\frac{1- F(u_n)}{v_n w(u_n)}\int_{ z^2/2+O(1/t_n^2)}^\IF
\frac{1- F(u_n+ s/w(u_n))}{1- F(u_n)} \psi_{\tau_n}(z^2/2(1+o(1)))\, ds\\
 &=  & (1+o(1))\frac{\rho }{ \sqrt{1- \rho^2}}\psi(z^2/2)h_n(1/\sqrt{t_n})\frac{1- F(u_n)}{t_n}\int_{ z^2/2}^\IF
\exp(-s) \, ds, \quad \ntoi.
\EQNY
Define next for $n\ge 1$ and $s\ge 0$
\def\gns{g_n(s)}
$$\gns:=  \frac{1 }{l_n(s)\sqrt{l_n^2(s)- 1}}.$$
We have for all $s\ge 0, n\ge 1$
\BQNY
\Abs{\gns-\wtiro}&< & (\frac{z}{\sqrt{t_n}}+ s/t_n)K, \quad \gns <  \wtiro K, \quad \wtiro:=\rho^2 /\sqrt{1- \rho^2},
\EQNY
with $K>1$ a positive constant. Next, for any $\zeta\ge 0$ we may write
\BQNY
\lefteqn{
\Abs{ \frac{v_n w(u_n)}{h_n(1/\sqrt{t_n})[1- F(u_n)]}J(l_n(\zeta) , v_n,\btrhoN )
- \wtiro \psi(\zeta)\exp(-\zeta)}} \\
 &=  &  \Abs{ \int_{\zeta }^\IF \frac{1- F(v_nl_n(s))}{1- F(u_n)}\frac{ \btrhoN(l_n(s)  )}{h_n(1/\sqrt{t_n})}
 g_n(s)\, ds - \tilde \rho \int_{ \zeta}^\IF \psi(\zeta)\exp(-s) \, ds}\\
 &\le & \int_{ \zeta }^\IF \Abs{ \frac{1- F(u_n+ s/w(u_n))}{1- F(u_n)}\frac{ \btrhoN(l_n(s))}{h_n(1/\sqrt{t_n})}
 g_n(s)\, ds - \wtiro \psi(\zeta)\exp(-s)} \, ds\\
 &\le & \int_{ \zeta}^\IF \Abs{ \frac{1- F(u_n+ s/w(u_n))}{1- F(u_n)}- \exp(-s)}
 \frac{ \btrhoN(l_n(s)  )}{h_n(1/\sqrt{t_n})} g_n(s)\, ds \\
 && +\int_{ \zeta }^\IF \exp(-s) \Abs{ \frac{ \btrhoN(l_n(s) )}{h_n(1/\sqrt{t_n})}
 g_n(s)- \wtiro \psi(\zeta)}\, ds \\
 &\le &  A(u_n) K \wtiro \psi(\zeta) \int_{ \zeta }^\IF B(s) \, ds +\wtiro \int_{ \zeta }^\IF \exp(-s) \Abs{ \frac{ \btrhoN(l_n(s) )}{h_n(1/\sqrt{t_n})} - \psi(\zeta)}\, ds \\
 && +  \int_{ \zeta }^\IF \exp(-s) \Abs{  g_n(s)- \wtiro }
 \frac{ \btrhoN(l_n(s)  )}{h_n(1/\sqrt{t_n})}\, ds \\
 &\le &   A(u_n) K \wtiro \psi(\zeta) \int_{ \zeta }^\IF B(s) \, ds
 +\wtiro a(t_n) b_n(\zeta) \exp(-\zeta)    +\psi(\zeta)\int_{ \zeta }^\IF \exp(-s) \Abs{  g_n(s)- \wtiro}\, ds \\
&=& R_n(u_n) .
\EQNY
Hence for all $n$ large since $\arhon= 1+ z^2/(2t_n)+ O(1/t_n^2)$ we may write
\BQNY
\Abs{ \frac{v_n w(u_n)}{h_n(1/\sqrt{t_n})(1- F(u_n))}J(\chi_n^{-1} \arhon, v_n,\btrhoN )
- \wtiro \psi(z^2/2)\exp(-z^2/2)}
&=& R_n(u_n).
\EQNY
Consequently, as $n\to \IF$
\BQNY
 J(\arhon,u_n, \wthN) &=&  h_n(1/\sqrt{t_n}) \frac{1- F(u_n)}{\sqrt{t_n}} \Biggl[[1- \Psi(z)]
 \int_{-\IF}^\IF\exp(-x^2/2) \psi(x^2/2) \, dx+ R_n(u_n)\Biggr]
\EQNY
and
\BQNY
J(\chi_n^{-1} \arhon, v_n,\btrhoN )
 &=& \frac{1}{\sqrt{t_n}}   h_n(1/\sqrt{t_n})\frac{1- F(u_n)}{\sqrt{t_n}}
 \Biggl[\frac{\rho }{ \sqrt{1- \rho^2}}\psi(z^2/2)\exp(-z^2/2) +R_n(u_n)\Biggr]
\EQNY
implying
\BQNY
\pk{Y_n> v_n \lvert X_n> u_n}&=& 1- \Psi(z) + \frac{1}{\sqrt{t_n}}\frac{\rho}{\sqrt{1- \rho^2} }\frac{ \exp(-z^2/2) \psi(z^2/2)}
{  \int_{-\IF}^\IF\exp(-x^2/2) \psi(x^2/2) \, dx}+ R_n(u_n), \quad \ntoi.
\EQNY
The proof for $z\le 0$ follows with similar calculations, hence the result. \QED

{\bf Acknowledgement:} I would like to thank two referees of Extremes for careful reading, several suggestions and corrections
received in April 2008.

\end{document}